\RequirePackage{etoolbox}
\csdef{input@path}{{style/}{graphics/}}
\documentclass[aos,MSNbibl,nameyear,dvips]{arximspdf}
\usepackage{graphicx}
\doi{10.1214/14-AOS1287}
\volume{43}
\issue{2}
\pubyear{2015}
\firstpage{501}
\lastpage{518}
\docsubty{FLA}

\makeatletter

\def\cal{\mathcal}

\def\R{{\mathbb{R}}}
\def\E{{\mathbb{E}}}
\def\I{{\mathbb{I}}}
\def\proj{\mathrm{pr}^{\bot}}

\def\argmax{\operatorname{argmax}}
\def\argmin{\operatorname{argmin}}

\newtheorem{theorem}{Theorem}
\newproclaim{remark}{Remark}
\newtheorem{lemma}{Lemma}
\newproclaim{definition}{Definition}
\newtheorem{corol}{Corollary}
\newproclaim{example}{Example}
\newtheorem{proposition}{Proposition}

\newproclaim{assumption}{Assumption}
\makeatother

\begin{document}
\begin{frontmatter}

\title{Universally optimal designs for two interference~models}
\runtitle{Design for interference model}

\begin{aug}
\author[A]{\fnms{Wei}~\snm{Zheng}\corref{}\ead
[label=e1]{weizheng@iupui.edu}}
\runauthor{W. Zheng}
\affiliation{Indiana University--Purdue University Indianapolis}
\address[A]{Department of Mathematical Sciences\\
Indiana University--Purdue University Indianapolis\\
Indianapolis, Indiana 46202-3216\\
USA\\
\printead{e1}}
\end{aug}

%
\received{\smonth{5} \syear{2014}}
%
\revised{\smonth{9} \syear{2014}}

%
\begin{abstract}
A systematic study is carried out regarding universally optimal designs
under the interference model, previously investigated by
Kunert and Martin
[\textit{Ann. Statist.} \textbf{28} (2000) 1728--1742]
and
Kunert and Mersmann
[\textit{J. Statist. Plann. Inference} \textbf{141} (2011) 1623--1632].
Parallel results are also
provided for the undirectional interference model, where the left and
right neighbor effects are equal. It is further shown that the
efficiency of any design under the latter model is at least its
efficiency under the former model. Designs universally optimal for both
models are also identified. Most importantly, this paper provides
Kushner's type linear equations system as a necessary and sufficient
condition for a design to be universally optimal. This result is novel
for models with at least two sets of treatment-related nuisance
parameters, which are left and right neighbor effects here. It sheds
light on other models in deriving asymmetric optimal or efficient designs.
\end{abstract}

%
\end{frontmatter}

\section{Introduction}
One issue with the application of block designs in agricultural field
trials is that a treatment assigned to a particular plot typically has
effects on the neighboring plots besides the effect on its own plot.
See \citet{Ree67}, \citet{DraGut80}, \citet{Kem82}, \citet{BesKem86}, \citet{Lan90}, \citet{Gil93},
Goldringer, Brabant and Kempton (\citeyear{GolBraKem94}),
Clarke, Baker and DePauw (\citeyear{ClaBakDeP00}), \citet{Davetal01} and \citet{Conetal08}
for examples in various backgrounds. Interference models have been
suggested for the analysis of data in order to avoid systematic bias
caused by these neighbor effects. Various designs have been proposed by
\citet{Gil93}, \citet{Dru99},
Filipiak and Markiewicz (\citeyear{FilMar03,FilMar05,FilMar07}), \citet{BaiDru04}, Ai, Ge and Chan (\citeyear{AiGeCha07}), Ai, Yu and He (\citeyear{AiYuHe09}),
Druilhet and Tinssonb (\citeyear{DruTin12}) and \citet{Fil12} among others. All of
them considered circular designs, where each block has a guard plot at
each end so that each plot within the block has two neighbors.

To study noncircular designs, \citet{KunMar00} investigated the
case when the block size, say $k$, is $3$ or $4$, which is extended by
\citet{KunMer11} to $t\geq k\geq5$, where $t$ is the number
of treatments. Both of them restricted to the subclass of pseudo
symmetric designs and the assumption that the within-block covariance
matrix is proportional to the identity matrix. This paper provides a
unified framework for deriving optimal pseudo symmetric designs for an
arbitrary covariance matrix as well as the general setup of $k\geq3$
and $t\geq2$. Most importantly, the Kushner's type linear equations
system is developed as a necessary and sufficient condition for any
design to be universally optimal, which is a powerful device for
deriving asymmetric designs. Moreover, a new approach of finding the
optimal sequences are proposed. These results are novel for models with
at least two sets of treatment-related nuisance parameters, which are
left and right neighbor effects here. They shed light on other similar
or more complicated models such as the one in Afsarinejad and Hedayat
(\citeyear{AfHe02}) and Kunert and Stufken (\citeyear{KuSt02}) for the
study of crossover
designs. Here, parallel results are also provided for the undirectional
interference model where the left and right neighbor effects are equal.
It is further established that the efficiency of any given design under
the latter model is not less than the one under the former model, for
the purpose of estimating the direct treatment effects.

Throughout the paper, we consider designs in $\Omega_{n,k,t}$, the set
of all possible block designs with $n$ blocks of size $k$ and $t$
treatments. The response, denoted as $y_{dij}$, observed from the $j$th
plot of block $i$ is modeled as
%
\begin{equation}
\label{eqn:729} Y_{dij}=\mu+\beta_i+\tau_{d(i,j)}+
\lambda_{d(i,j-1)}+\rho _{d(i,j+1)}+\varepsilon_{ij},
\end{equation}
where $\E\varepsilon_{ij}=0$. The subscript $d(i,j)$ denotes the
treatment assigned in the $j$th plot of block $i$ by the design $d\dvtx \{
1,2,\ldots,n\}\times\{1,2,\ldots,k\}\rightarrow
\{1,2,\ldots,t\}$. Furthermore, $\mu$ is the general mean, $\beta_i$ is
the $i$th block effect, $\tau_{d(i,j)}$ is the direct treatment effect
of treatment $d(i,j)$, $\lambda_{d(i,j-1)}$ is the neighbor effect of
treatment $d(i,j-1)$ from the left neighbor, and $\rho_{d(i,j+1)}$ is
the neighbor effect of treatment $d(i,j+1)$ from the right neighbor.
One major objective of design theorists is to find optimal or efficient
designs for estimating the direct treatment effects in the model.

If $Y_d$ is the vector of responses organized block by block, model
(\ref{eqn:729}) is written in a matrix form of
%
\begin{equation}
\label{eqn:728} Y_d =1_{nk}\mu+U\beta+T_d
\tau+L_d\lambda+R_d\rho+\varepsilon,
\end{equation}
where $\beta=(\beta_1,\ldots,\beta_n)'$, $\tau=(\tau_1,\ldots,\tau_t)'$,
$\lambda=(\lambda_1,\ldots,\lambda_t)'$ and $\rho=(\rho_1,\ldots,\rho_t)'$.
The notation $'$ means the transpose of a vector or a matrix. Here, we
have $U=I_n\otimes1_k$ with $\otimes$ as the Kronecker product, and
$1_k$ represents a vector of ones with length $k$. Also, $T_d$, $L_d$
and $R_d$ represent the design matrices for the direct, left neighbor
and right neighbor effects, respectively. We assume there is no guard
plots, that is, $\lambda_{d(i,0)}=\rho_{d(i,k+1)}=0$. Then we have
$L_d=(I_n\otimes H) T_d$ and $R_d=(I_n\otimes H') T_d$, where
$H(i,j)=\I
_{i=j+1}$ with the indicator function $\I$.

Here, we merely assume $\operatorname{Var}(\varepsilon)=I_n\otimes\Sigma$, with
$\Sigma$ being an arbitrary $k\times k$ positive definite symmetric
matrix. Given a matrix, say $G$, we define the projection
$\proj G=I-G(G'G)^-G'$. The information matrix for the direct treatment
effect $\tau$ is
%
\begin{equation}
\label{eqn:912} C_d=T_d'V'
\proj(VU|VL_d|VR_d)VT_d,
\end{equation}
where $V$ is the matrix such that $V^2=I_n\otimes\Sigma^{-1}$. By
direct calculations, we have
\begin{eqnarray*}
C_d&=&E_{d00}-E_{d01}E_{d11}^-E_{d10},
\\
E_{d00}&=&C_{d00},
\\
E_{d10}'&=& E_{d01}=\pmatrix{C_{d01} & C_{d02}},
\\
E_{d11}&=& \pmatrix{ C_{d11} &
C_{d12}
\cr
C_{d21} & C_{d22}},
\end{eqnarray*}
where $C_{dij}=G_i'(I_n\otimes\tilde{B})G_j,0\leq i,j\leq2$ with
$G_0=T_d$, $G_1=L_d$, $G_2=R_d$ and $\tilde{B}=\Sigma^{-1}-\Sigma
^{-1}J_k\Sigma^{-1}/1_k'\Sigma^{-1}1_k$ with $J_k=1_k1_k'$. It is
obvious that $C_{dij}=C_{dji}'$. For the special case of $\Sigma=I_k$,
we have the simplification of $\tilde{B}=I_k-k^{-1}J_k=\proj(1_k)$, and
the latter is denoted by $B_k$. \citet{Kus97} pointed out that when
$\Sigma$ is of type-$H$, that is, $aI_k+b1_k'+1_kb'$ with $a\in\R$ and
$b\in\R^k$, we have
%
\begin{equation}
\label{eqn:05233} \tilde{B}=B_k/a.
\end{equation}
Hence, the choices of designs agree with that for $\Sigma=I_k$. This
special case will be particularly dealt with in Section~\ref{sec5}. We
allow $\Sigma$ to be an arbitrary covariance matrix throughout the rest
of the paper.

Note that a design in $\Omega_{n,k,t}$ could be considered as a result
of selecting $n$ elements from the set, ${\cal S}$, of all possible
$t^k$ block sequences with replacement. For sequence $s\in{\cal S}$,
define the sequence proportion $p_s=n_s/n$, where $n_s$ is the number
of replications of $s$ in the design. A design is determined by
$n_s,s\in{\cal S}$, which is in turn determined by the \textit{measure}
$\xi=(p_s,s\in{\cal S})$ for any fixed $n$.

For $0\leq i,j\leq2$, define $C_{sij}$ to be $C_{dij}$ when the design
consists of the single sequence $s$, and let $C_{\xi ij}=\sum_{s\in
{\cal S}}p_sC_{sij}$. Then we have $C_{dij}=nC_{\xi ij},0\leq i,j\leq
2$. Similarly, $E_{dij}=n\sum_{s\in{\cal S}}p_sE_{sij}=nE_{\xi
ij},0\leq i,j\leq1$. Note that $C_d$ is a Schur's complement of
$A_d=(E_{dij})_{0\leq i,j\leq1}$, for which we also have $A_d=n\sum_{s\in{\cal S}}p_sA_s=nA_{\xi}$. It is obvious that $C_d=nC_{\xi}$,
where $C_{\xi}=E_{\xi00}-\break E_{\xi01}E_{\xi11}^-E_{\xi10}$. In
approximate design theory, we try to find the optimal measure $\xi$
among the set ${\cal P}=\{(p_s,s\in{\cal S})|\sum_{s\in{\cal
S}}p_s=1,p_s\geq0\}$ to maximize $\Phi(C_{\xi})$ for a given function
$\Phi$ satisfying the following three conditions [\citet{Kie75}]:
\begin{longlist}[(C.1)]
\item[(C.1)] $\Phi$ is concave.
\item[(C.2)] $\Phi(S'CS)=\Phi(C)$ for any permutation matrix $S$.
\item[(C.3)] $\Phi(bC)$ is nondecreasing in the scalar $b>0$.
\end{longlist}
A measure $\xi$ which achieves the maximum of $\Phi(C_{\xi})$ among
${\cal P}$ for any $\Phi$ satisfying (C.1)--(C.3) is said to be \textit{universally optimal}.
Such measure is optimal under criteria of $A$, $D$,
$E$, $T$, etc.

The rest of the paper is organized as follows. Section~\ref{sec2}
provides some preliminary results as well as a necessary and sufficient
condition for a pseudo symmetric measure to be universally optimal
among ${\cal P}$. The latter is critical for deriving the optimal
sequence proportions through an algorithm. Section~\ref{sec3} provides
a linear equations system of $p_s,s\in{\cal S}$, as a necessary and
sufficient condition for a measure to be universally optimal. Section~\ref{sec4} provides similar results for the model with $\lambda=\rho$.
Further, it is shown that the efficiency of any design under the latter
model would be at least its efficiency under model (\ref{eqn:728}).
Also, an alternative approach is given to derive the optimal sequence
proportions. Section~\ref{sec5} derives theoretical results regarding
feasible sequences when $\Sigma$ is of type-$H$. Section~\ref{sec6}
provides some examples of optimal or efficient designs for various
combinations of $k,t,n$ and $\Sigma$.

\section{Pseudo symmetric measure}\label{sec2}
Let ${\cal G}$ be the set of all $t!$ permutations on symbols $\{
1,2,\ldots,t\}$. For permutation $\sigma\in{\cal G}$ and sequence
$s=(t_1\cdots t_k)$ with $1\leq t_i\leq t$ and $1\leq i\leq k$, we
define $\sigma s=(\sigma(t_1)\cdots\sigma(t_k))$. For measure $\xi
=(p_s,s\in{\cal S})$, we define $\sigma\xi=(p_{\sigma^{-1}s},s\in
{\cal S})$. A measure is said to be \textit{symmetric} if $\sigma\xi
=\xi$
for all $\sigma\in{\cal G}$. For sequence $s$, denote by $\langle
s\rangle=\{\sigma s\dvtx \sigma\in{\cal G}\}$ the \textit{symmetric block}
generated by $s$. Such symmetric blocks are also called equivalent
classes by \citet{Kus97}, due to the fact that symmetric blocks
generated by two different sequences are either identical or mutually
disjoint. Now let $m$ be the total number of distinct symmetric blocks
which partition ${\cal S}$. Without loss of generality, suppose these
$m$ symmetric blocks are generated by sequences $s_i$, $1\leq i\leq m$.
Then we have ${\cal S}=\bigcup^m_{i=1}\langle s_i\rangle$. For a symmetric
measure, we have
%
\begin{equation}
\label{eqn:0325} p_s=p_{\langle s_i\rangle}/\bigl|\langle s_i
\rangle\bigr|\qquad \mbox{for } s\in\langle s_i\rangle, 1\leq i\leq m,
\end{equation}
where $p_{\langle s_i\rangle}=\sum_{s\in\langle s_i\rangle}p_s$ and
$|\langle s_i\rangle|$ is the cardinality of $\langle s_i\rangle$. The
linearity of $A_d$, conditions (C.1)--(C.3) and properties of Schur's
complement together yield the following lemma.

\begin{lemma}\label{lemma:0223}
For any measure, say $\xi$, there exists a symmetric measure, say $\xi
^*$, such that $\Phi(C_{\xi})\leq\Phi(C_{\xi^*})$ for any $\Phi$
satisfying \textup{(C.1)--(C.3)}.
\end{lemma}

Define a measure to be \textit{pseudo symmetric} if $C_{\xi ij},0\leq
i,j\leq2$ are all completely symmetric. It is easy to verify that a
symmetric measure is also pseudo symmetric. The difference is that
(\ref
{eqn:0325}) does not has to hold for a general pseudo symmetric
measure. Lemma~\ref{lemma:0223} indicates that an optimal measure in
the subclass of (pseudo) symmetric measures is automatically optimal
among ${\cal P}$. For a pseudo symmetric measure, we have $C_{\xi
ij}=c_{\xi ij}B_t/(t-1)+(1_t'C_{\xi ij}1_t)J_t/t^2$, $0\leq i,j\leq2$,
where $c_{\xi ij}=\operatorname{tr}(B_t C_{\xi ij}B_t)$. Hence $E_{\xi11}=Q_{\xi
}\otimes B_t/(t-1)+\tilde{Q}_{\xi}\otimes J_t/t^2$, where $Q_{\xi
}=(c_{\xi ij})_{1\leq i,j\leq2}$ and $\tilde{Q}_{\xi}=(1_t'C_{\xi
ij}1_t)_{1\leq i,j\leq2}$. Now we show that both $Q_{\xi}$ and
$\tilde
{Q}_{\xi}$ are positive definite for any measure, and hence $E_{\xi
11}$ is positive definite for any pseudo symmetric measure. The latter
is the key to prove Theorem~\ref{thm:0325}.

\begin{lemma}\label{lemma323}
$Q_{\xi}$ is positive definite for any measure $\xi$.
\end{lemma}

\begin{pf}It is sufficient to show the nonsingularity of $Q_s$ for
all $s\in{\cal S}$. Suppose $Q_s$ is singular, there exists a nonzero
vector $x=(x_1,x_2)'$ such that
\begin{eqnarray*}
0&=&x'Q_sx=\sum_{i=1}^2
\sum_{j=1}^2 x_ix_jc_{s ij}
\\
&=&\operatorname{tr} \Biggl(\sum_{i=1}^2\sum
_{j=1}^2 x_ix_jB_tC_{s ij}B_t
\Biggr).
\end{eqnarray*}
Since $\sum_{i=1}^2\sum_{j=1}^2 x_ix_jB_tC_{s ij}B_t$ is a nonnegative
definite matrix, we have
\begin{eqnarray*}
0&=&\sum_{i=1}^2\sum
_{j=1}^2 x_ix_jB_tC_{s ij}B_t
\\
&=&B_t(x_1L_s+x_2R_s)'
\tilde{B}(x_1L_s+x_2R_s)B_t,
\end{eqnarray*}
which in turn yields
%
\begin{equation}
\label{eqn:323} 0=\tilde{B}(x_1L_s+x_2R_s)B_t.
\end{equation}
Equation (\ref{eqn:323}) is only possible when each column of
$M=(x_1L_s+x_2R_s)B_t$ consists of identical entries, that is, the rows
of $M$ are identical. In the sequel, we investigate the possibility of
(\ref{eqn:323}) for sequence $s=(t_1\cdots t_k)$. Define $e_i$ to be a
zero--one vector of length $t$ with only its $i$th entry as one, then
the first, second and last rows of $M$ are given by
$x_2(e_{t_2}-1_t/t)'$, $x_1(e_{t_1}-1_t/t)'+x_2(e_{t_3}-1_t/t)'$ and
$x_1(e_{t_{k-1}}-1_t/t)'$, respectively. Now we continue the discussion
in the following four cases. (i) If $x_1=x_2$, the equality of the
first two rows of $M$ indicates $e_{t_1}+e_{t_3}-e_{t_2}=1_t/t$, which
is impossible since the left-hand side is a vector of integers and the
right-hand side is a vector of fractional numbers. (ii) If $x_1\neq
x_2$ and $t_2=t_{k-1}$, the first and the last rows of $M$ cannot be
the same. (iii) If $x_1\neq x_2$, $t_2\neq t_{k-1}$ and $t=2$, the
equality of the first and the last rows of $M$ necessities $x_1+x_2=0$,
which together with the equality of the first two rows of $M$ indicates
$e_{t_1}+e_{t_2}-e_{t_3}=1_t/t$, which is again impossible. (iv) If
$x_1\neq x_2$, $t_2\neq t_{k-1}$ and $t\geq3$, by looking at the
$t_2$th and $t_{k-1}$th entries of the first and last rows of $M$,
(\ref
{eqn:323}) necessities $x_2(1-1/t)=-x_1/t$ and $x_1(1-1/t)=-x_2/t$
which is impossible by simple algebra.
\end{pf}
%

\begin{lemma}\label{lemma325}
$\tilde{Q}_{\xi}$ is positive definite for any measure $\xi$.
\end{lemma}
\begin{pf}Since $\tilde{B}$ has column and row sums as zero. We have
%
\begin{equation}
\tilde{Q}_{\xi}=\pmatrix{
\tilde{B}(1,1) & \tilde{B}(1,k)
\cr
\tilde{B}(k,1) & \tilde{B}(k,k)},
\end{equation}
where $\tilde{B}(i,j)$ means the $(i,j)$th entry in $\tilde{B}$. For
vector $x=(x_1,x_2)'\in\R^2$, define $w=(x_1,0,\ldots,0,x_2)'\in\R^k$.
For any nonzero $x$, we have
%
\begin{equation}
x'\tilde{Q}_{\xi}x=w'\tilde{B}w>0,
\end{equation}
in view of the fact that $\tilde{B}1_k=0$, $\tilde{B}\geq0$ and the
rank of $\tilde{B}$ is $k-1$. Hence, the lemma is concluded.
\end{pf}

\begin{lemma}\label{lemma:0205}
For a pseudo symmetric measure, say $\xi$, we have $C_{\xi}=q_{\xi
}^*B_t/\break (t-1)$, where
%
\begin{equation}
\label{eqn:1109} q_{\xi}^*=c_{\xi00}-\ell_{\xi}'Q_{\xi}^{-1}
\ell_{\xi},
\end{equation}
with $\ell_{\xi}=(c_{\xi01}, c_{\xi02})'$.
\end{lemma}

\begin{remark}
In proving Lemma~\ref{lemma:0205}, we used the equations $1_t'C_{\xi
0j}=0$, $0\leq j\leq2$. Note that $nq_{\xi}^*$ is the $q_d^*$ as
defined in \citet{KunMar00}. Lemma~\ref{lemma323} shows that
only case (i) of the four cases proposed by them is possible. Hence the
generalized inverse $Q_{\xi}^-$ in \citet{KunMar00} is now
replaced by $Q_{\xi}^{-1}$ in~(\ref{eqn:1109}).
\end{remark}

By applying Lemmas \ref{lemma:0223} and \ref{lemma:0205}, we derive the
following proposition.

\begin{proposition}\label{prop:323}
Let $y^*=\max_{\xi\in{\cal P}}q_{\xi}^*$. A measure $\xi\in{\cal
P}$ is universally optimal \textup{(i)} if it is a pseudo symmetric measure with
$q_{\xi}^*=y^*$, \textup{(ii)} if and only if $C_{\xi}=y^*B_t/(t-1)$.
\end{proposition}

Let $R_s=(c_{sij})_{0\leq i,j\leq2}$ and $R_{\xi}=\sum_{s\in{\cal
S}}p_sR_s$. By Lemma~\ref{lemma323} we have $q_{\xi}^*=\operatorname{det}(R_{\xi
})/\operatorname{det}(Q_{\xi})$, where $\operatorname{det}(\cdot)$ means the determinant of a square
matrix. For measure $\xi=(p_s, s\in{\cal S})$, we call the set ${\cal
V}_{\xi}=\{s\dvtx p_s>0,s\in{\cal S}\}$ the \textit{support} of~$\xi$. One
can identify universally optimal pseudo symmetric measures based on the
following theorem. See \citet{Zhe13N2} for an algorithm based on a
similar theorem.

\begin{theorem}\label{thm:323}
A pseudo symmetric measure, say $\xi$, is universally optimal if and
only if $\operatorname{det}(R_{\xi})>0$ and
%
\begin{equation}
\label{eqn:3232} \max_{s\in{\cal S}} \bigl[\operatorname{tr}\bigl(R_sR_{\xi}^{-1}
\bigr)-\operatorname{tr}\bigl(Q_sQ_{\xi
}^{-1}\bigr) \bigr]=1.
\end{equation}
Moreover, each sequence in ${\cal V}_{\xi}$ reaches the maximum in
(\ref
{eqn:3232}).
\end{theorem}

\begin{pf}
If $\operatorname{det}(R_{\xi})=0$, we have $q_{\xi}^*=0$, which means that such
design has no information regarding $\tau$, and hence can be readily
excluded from the consideration. In the sequel, we restrict the
discussion to the case of $\operatorname{det}(R_{\xi})>0$.

By Lemmas \ref{lemma:0223}, \ref{lemma323} and \ref{lemma:0205}, a
pseudo symmetric measure, say $\xi$, is universally optimal if and only
if it achieves the maximum of $\varphi(\xi)=\operatorname{log}(\operatorname{det}(R_{\xi
})/\operatorname{det}(Q_{\xi
}))$, which is equivalent to
%
\begin{equation}
\label{eqn:8014} \lim_{\delta\rightarrow0} \frac{\varphi[(1-\delta)\xi+\delta\xi
_0]-\varphi(\xi)}{\delta}\leq0,
\end{equation}
for any measure $\xi_0\in{\cal P}$. It is well known that
%
\begin{equation}
\label{eqn:8013} \lim_{\delta\rightarrow0} \frac{\operatorname{log}(\operatorname{det}(R_{(1-\delta)\xi+\delta
\xi
_0}))-\operatorname{log}(\operatorname{det}(R_{\xi}))}{\delta}=\operatorname{tr}
\bigl(R_{\xi_0}R_{\xi}^{-1}\bigr)-3.
\end{equation}
The same result holds for $Q(\xi)$ except that $3$ should be replaced
by $2$. By applying~(\ref{eqn:8013}) to (\ref{eqn:8014}), we have
%
\begin{equation}
\label{eqn:8015} \operatorname{tr}\bigl(R_{\xi_0}R_{\xi}^{-1}\bigr)-\operatorname{tr}
\bigl(Q_{\xi_0}Q_{\xi}^{-1}\bigr)\leq 1.
\end{equation}
In (\ref{eqn:8015}), by setting $\xi_0$ to be a degenerated measure
which puts all its mass on a single sequence, we derive
\[
\max_{s\in{\cal S}} \bigl(\operatorname{tr}\bigl(R_sR_{\xi}^{-1}
\bigr)-\operatorname{tr}\bigl(Q_sQ_{\xi
}^{-1}\bigr) \bigr)\leq 1.
\]
By taking $\xi_0=\xi$, we have the equal sign for (\ref{eqn:8015}).
Also observe that conditioning on fixed $\xi$, the left-hand side of
(\ref{eqn:8015}) is a linear function of the proportions in $\xi_0$.
Thus, we have
\[
\max_{s\in{\cal S}} \bigl(\operatorname{tr}\bigl(R_sR_{\xi}^{-1}
\bigr)-\operatorname{tr}\bigl(Q_sQ_{\xi
}^{-1}\bigr) \bigr)\geq 1.
\]
Hence, the theorem follows.
\end{pf}

\section{The linear equations system: A necessary and sufficient
condition for universal optimality}\label{sec3}
For sequence $s$ and vector $x\in\R^2$, define the quadratic function
$q_s(x)=c_{s00}+2\ell_s'x+x'Q_sx$. For measure $\xi=(p_s,s\in{\cal
S})$, define $q_{\xi}(x)=\sum_{s\in{\cal S}}p_sq_s(x)=c_{\xi
00}+2\ell
_{\xi}'x+x'Q_{\xi}x$. One can verify that $q_{\xi}^*=\break \min_{x\in\R
^2}q_{\xi}(x)$. Since $q_s(x)$ is strictly convex for all $s\in{\cal
S}$ in view of Lemma~\ref{lemma323}, thus $r(x):=\max_{s\in{\cal
S}}q_s(x)$ is also strictly convex. Let $x^*$ be the unique point in
$\R
^2$ which achieves minimum of $r(x)$ and define ${\cal T}=\{
s\dvtx q_s(x^*)=r(x^*), s\in{\cal S}\}$. Recall $y^*=\max_{\xi\in{\cal
P}}q_{\xi}^*$ and ${\cal V}_{\xi}=\{s\dvtx p_s>0,s\in{\cal S}\}$, now we
derive Theorem~\ref{thm:325} below which is important for proving
Theorem~\ref{thm:0325} and results in Section~\ref{sec4}.

\begin{theorem}\label{thm:325}
\textup{(i)} $y^*=r(x^*)$. \textup{(ii)} $q^*_{\xi}=y^*$ implies $x^*=-Q_{\xi}^{-1}\ell
_{\xi}$. \textup{(iii)} $q^*_{\xi}=y^*$ implies ${\cal V}_{\xi}\subset{\cal T}$.
\end{theorem}

\begin{pf} First, we have
\[
y^*=\max_{\xi\in{\cal P}}\min_{x\in\R
^2}q_{\xi}(x)\leq \min_{x\in\R^2}\max_{\xi\in{\cal P}}q_{\xi}(x)
=\min_{x\in\R^2}\max_{s\in{\cal S}}q_s(x)=r(x^*).
\]
Then (i) is
proved if we can show $y^*\geq r(x^*)$. To see the latter, define
${\cal T}_0=\{s\dvtx q_s(x^*)=r(x^*), s\in\{s_1\cdots s_m\}\}$. $(1)$ If
${\cal T}_0$ contains a single sequence, say $s_1$, let $\xi_0$ be the
measure with $p_{\langle s_1\rangle}=1$, then we have $\min_{x\in\R
^2}q_{\xi_0}(x)=r(x^*)$. Hence, $y^*\geq r(x^*)$. $(2)$ If ${\cal T}_0$
contains more than one sequences, let $\nabla q_s(x^*)$ be the gradient
of $q_s(x)$ evaluated at point $x=x^*$ and define $\Xi$ to be the
convex hull of $\{\nabla q_s(x^*)\dvtx s\in{\cal T}_0\}$. We claim $0\in
\Xi
$, since otherwise we could find a vector $z\in\R^2$ so that
$z'\nabla
q_s(x^*)<0$ for all $s\in\{\nabla q_s(x^*)\dvtx s\in{\cal T}_0\}$,\vadjust{\goodbreak} which
would indicate that $x^*$ is not the minimum point of $r(x)$, and hence
the contradiction is reached. Note that $0\in\Xi$ indicates there
exists a measure, say $\xi_0$, such that $q_{\xi_0}(x^*)=r(x^*)$ and
$\nabla q_{\xi_0}(x^*)=0$, which yields $\min_{x\in\R^2}q_{\xi
_0}(x)=r(x^*)$ and hence $y^*\geq r(x^*)$. (i) is thus proved.

Observe that the minimum of $q_{\xi}(x)$ is achieved at the unique
point $x=-Q_{\xi}^{-1}\ell_{\xi}:=\tilde{x}$. If $\tilde{x}\neq x^*$,
we have $y^*=r(x^*)\geq q_{\xi}(x^*)>q_{\xi}(\tilde{x})=q^*_{\xi}$ and
hence the contradiction is reached. (ii) is thus concluded.

For (iii), if there is a sequence, say $s$, with $s\in{\cal V}_{\xi
}$ and $s\notin{\cal T}$, we have $y^*>q_{\xi}(x^*)\geq q^*_{\xi}$,
and hence the contradiction is reached.
\end{pf}

\begin{theorem}\label{thm:0325}
A measure $\xi=(p_s,s\in{\cal S})$ is universally optimal among
${\cal
P}$ if and only if
%
\begin{eqnarray}
\sum_{s\in{\cal T}}p_s\bigl[E_{s00}+E_{s 01}
\bigl(x^*\otimes B_t\bigr)\bigr]&=&y^*B_t/(t-1),\label{eqn:03252}
\\
\sum_{s\in{\cal T}}p_s\bigl[E_{s 10}+E_{s 11}
\bigl(x^*\otimes B_t\bigr)\bigr]&=&0,\label
{eqn:03253}
\\
\sum_{s\notin{\cal T}}p_s&=&0.\label{eqn:03254}
\end{eqnarray}
\end{theorem}

\begin{pf} Note that (\ref{eqn:03252})--(\ref{eqn:03254}) is
equivalent to
%
\begin{eqnarray}
E_{\xi00}+E_{\xi01}\bigl(x^*\otimes B_t
\bigr)&=&y^*B_t/(t-1),\label{eqn:3255}
\\
E_{\xi10}+E_{\xi11}\bigl(x^*\otimes B_t
\bigr)&=&0,\label{eqn:3254}
\\
\sum_{s\in{\cal T}}p_s&=&1.\label{eqn:03255}
\end{eqnarray}

{Necessity}. By Proposition~\ref{prop:323}, there exists a symmetric
measure, say $\xi_1$, which is universally optimal. Further, we have
$C_{\xi}=C_{\xi_1}=y^*B_t/(t-1)$. Define $\xi_2=(\xi+\xi_1)/2$.
Then we
have $A_{\xi_2}=(A_{\xi}+A_{\xi_1})/2$, which indicates $C_{\xi
_2}\geq
(C_{\xi}+C_{\xi_1})/2=y^*B_t/(t-1)$. The latter combined with
Proposition~\ref{prop:323} yields $C_{\xi_2}=y^*B_t/(t-1)$. Hence, by
similar arguments as in \citet{Kus97}, we have
%
\begin{eqnarray}
E_{\xi11}\bigl(E_{\xi11}^+E_{\xi10}-E_{\xi_211}^+E_{\xi
_210}
\bigr)&=&0,\label
{eqn:325}
\\
E_{\xi_111}\bigl(E_{\xi_111}^+E_{\xi_110}-E_{\xi_211}^+E_{\xi
_210}
\bigr)&=&0,\label{eqn:3252}
\end{eqnarray}
where $^+$ means the Moore--Penrose generalized inverse. Since $\xi_1$
is a symmetric measure, we have $E_{\xi_111}=Q_{\xi_1}\otimes
B_t/(t-1)+\tilde{Q}_{\xi_1}\otimes J_t/t^2$. By Lemmas \ref{lemma323},
\ref{lemma325} and the orthogonality between $B_t$ and $J_t$, we obtain
$\operatorname{det}(E_{\xi_111})=\operatorname{det}(Q_{\xi})^{t-1}\operatorname{det}(\tilde{Q}_{\xi
})/[(t-1)^{2t-2}t^3]>0$. Applying the latter to (\ref{eqn:3252}) yields
%
\begin{eqnarray}\label{eqn:3253}
E_{\xi_211}^+E_{\xi_210}&=&E_{\xi_111}^+E_{\xi_110}\nonumber
\\
&=&Q_{\xi_1}^{-1}\ell_{\xi_1}\otimes B_t
\\
&=&-x^*\otimes B_t,
\nonumber
\end{eqnarray}
in view of Theorem~\ref{thm:325}(ii). Now (\ref{eqn:3254}) is
derived from (\ref{eqn:325}) and (\ref{eqn:3253}). By (\ref{eqn:3254}),
we have
%
\begin{eqnarray}
y^*B_t/(t-1)&=&C_{\xi}=E_{\xi00}-E_{\xi01}E_{\xi11}^-E_{\xi
10}\label
{eqn:3256}
\\
&=&E_{\xi00}+E_{\xi01}E_{\xi11}^-E_{\xi11}\bigl(x^*
\otimes B_t\bigr)\label
{eqn:3257}
\\
&=&E_{\xi00}+E_{\xi01}\bigl(x^*\otimes B_t
\bigr),\label{eqn:3258}
\end{eqnarray}
which is essentially (\ref{eqn:3255}).

By (5.2) of \citet{Kus97}, we have $C_{\xi}\leq H'A_{\xi}H$ for any
$3t\times t$ matrix $H$. Set $H=(x_0,x_1,x_3)'\otimes B_t$ with
$x_0\equiv1$, we have
%
\begin{equation}
\label{eqn:03256} C_{\xi}\leq\sum^2_{i=0}
\sum^2_{j=0}x_ix_jB_tC_{\xi ij}B_t.
\end{equation}
By taking the trace of both sides of (\ref{eqn:03256}), we have
\begin{eqnarray*}
\operatorname{tr}(C_{\xi})&\leq&\sum^2_{i=0}
\sum^2_{j=0}x_ix_jc_{\xi ij}
\\
&=&q_{\xi}(x),
\end{eqnarray*}
for $x=(x_1,x_2)'$. Now set $x=-Q_{\xi}^{-1}\ell_{\xi}$, we have
$\operatorname{tr}(C_{\xi})\leq q^*_{\xi}\leq y^*$. Note that $\operatorname{tr}(C_{\xi})=y^*$ in
view of Proposition~\ref{prop:323}(ii). As a result, we have
$q^*_{\xi
}=y^*$ and thus~(\ref{eqn:03255}) in view of Theorem~\ref{thm:325}(iii).

Sufficiency of (\ref{eqn:3255})--(\ref{eqn:03255}) is trivial in view
of (\ref{eqn:3256})--(\ref{eqn:3258}).
\end{pf}

\section{Undirectional interference model}\label{sec4}
In many occasions, it is reasonable to believe that the neighbor
effects of each treatment from the left and the right should be the
same, that is, $\lambda=\rho$. With this condition, model (\ref
{eqn:728}) reduces to
%
\begin{equation}
\label{eqn:0413} Y_d =1_{nk}\mu+U\beta+T_d
\tau+(L_d+R_d)\lambda+\varepsilon.
\end{equation}
The information matrix, $\tilde{C}_d$, for $\tau$ under model (\ref
{eqn:0413}) is given by
\begin{eqnarray*}
\tilde{C}_d&=&C_{d00}-\tilde{C}_{d01}
\tilde{C}_{d11}^-\tilde {C}_{d10},
\\
\tilde{C}_{d10}'&=&\tilde{C}_{d01}=T_d'(I_n
\otimes\tilde {B}) (L_d+R_d),
\\
\tilde{C}_{d11}&=&(L_d+R_d)'(I_n
\otimes\tilde{B}) (L_d+R_d).
\end{eqnarray*}
It is obvious that $\tilde{C}_d/n$ only depends on the measure $\xi
=(p_s,s\in{\cal S})$, and we denote such matrix by $\tilde{C}_{\xi}$.
Let $\tilde{q}_s(z)=q_s((z,z)')$ and $\tilde{r}(z)=\max_{s\in{\cal
S}}\tilde{q}_s(z)$ for $z\in\R$. Note that $\tilde{r}(z)$ is strictly
convex due to the strict convexity of $r(x)$, hence there is an unique
minimizer of $\tilde{r}(z)$ which is denoted by $z^*$ here. By
following similar arguments as in Sections~\ref{sec2} and \ref{sec3},
one can derive the following theorem for universally optimal measures
under model (\ref{eqn:0413}) in view of Lemma~\ref{lemma:0413}(ii).

\begin{theorem}\label{thm:0413}
Let $y_0=\tilde{r}(z^*)$ and ${\cal T}_0=\{s\in{\cal S}\dvtx \tilde
{q}_s(z^*)=y_0\}$. For measure $\xi=(p_s,s\in{\cal S})$, the following
three sets of conditions are equivalent.
\textup{(i)} $\xi$ is universally optimal. \textup{(ii)} $\tilde{C}_{\xi}=y_0B_t/(t-1)$.
\textup{(iii)}
%
\begin{eqnarray}
\sum_{s\in{\cal T}_0}p_s\bigl[C_{s00}+z^*
\tilde{C}_{s
01}B_t\bigr]&=&y_0B_t/(t-1),\label{eqn:04132}
\\
\sum_{s\in{\cal T}_0}p_s\bigl[
\tilde{C}_{s10}+z^*\tilde{C}_{s
11}B_t\bigr]&=&0,
\\
\sum_{s\in{\cal T}_0}p_s&=&1.\label{eqn:04133}
\end{eqnarray}
\end{theorem}

The following lemma is the key to build up the connections between the
two models as given by Theorem~\ref{thm:0424}.

\begin{lemma}\label{lemma:0413}
If $\Sigma$ is persymmetric, we have the following. \textup{(i)}
$x^*=(z^*,z^*)'$. \textup{(ii)}~$y^*=y_0$. \textup{(iii)} ${\cal T}={\cal T}_0$.
\end{lemma}

\begin{pf}
For sequence $s=(t_1t_2\cdots t_p)$, define its
\textit{dual}
sequence as $s'=(t_p,t_{p-1}\cdots t_1)$. First we claim that
%
\begin{eqnarray}
\ell_s&=&\Lambda_2 \ell_{s'},\label{eqn:04243}
\\
Q_s&=&\Lambda_2 Q_{s'} \Lambda_2,\label{eqn:04244}
\end{eqnarray}
where $\Delta_h=(\I_{i+j=h+1})_{1\leq i,j\leq h}$. Then the function
$r(x)$ is symmetric about the line $x_1=x_2$, where $x=(x_1,x_2)'$.
This indicates that the two components of $x^*\in\R^2$ are identical.
From this, (i) and (ii) follows immediately. (iii) follows
directly from~(i) and (ii) by definitions of ${\cal T}$ and ${\cal T}_0$.

To prove (\ref{eqn:04243}) and (\ref{eqn:04244}), it is sufficient to
show $L_s=\Delta_kR_{s'}$, $R_s=\Delta_kL_{s'}$ and $\Delta_k \tilde
{B}\Delta_k=\tilde{B}$. The first two equations are trivial. To see the
latter, note that the persymmetry (and hence the bisymmetry) of $\Sigma
$ indicates the bisymmetry of $\Sigma^{-1}$ in view of Laplace's
formula for calculating the matrix inverse. Hence, the sum of the $i$th
column (or row) of $\Sigma^{-1}$ is equal to the sum of its $(k+1-i)$th
column, which indicates the bisymmetry of $\Sigma^{-1}J_k\Sigma^{-1}$,
and hence the bisymmetry of $\tilde{B}$.
\end{pf}

\begin{remark}
There is a wide range of covariance matrices which are persymmetric.
Examples include the identity matrix, the completely symmetric matrix,
the $\operatorname{AR}(1)$ type covariance matrix and the one used in Section~\ref{sec6}. By Corollary~2.2 of \citet{Kus97}, Lemma~\ref{lemma:0413}
still holds if $\Sigma=\Sigma_0+\gamma1_k'+1_k\gamma'$ with $\Sigma_0$
being persymmetric. In fact, the lemma holds as long as $\tilde{B}$ is
persymmetric. When $\tilde{B}$ is not persymmetric, empirical evidence
indicates that we typically have $x^*\neq(z^*,z^*)'$ and $y^*<y_0$.
Even though we observe ${\cal T}={\cal T}_0$ very often, however, the
optimal proportions for sequences in the support would be different for
the two models.
\end{remark}

A measure $\xi=(p_s,s\in{\cal S})$ is said to be \textit{dual} if
$p_{\langle s\rangle}=p_{\langle s'\rangle}, s\in{\cal S}$, where $s'$
is the dual sequence of $s$ as defined in the proof of Lemma \ref{lemma:0413}.

\begin{theorem}\label{thm:0424}
If $\Sigma$ is persymmetric, we have the following. \textup{(i)} For any
measure, its universal optimality under model (\ref{eqn:728}) implies
its universal optimality under model (\ref{eqn:0413}). \textup{(ii)} For a
pseudo symmetric dual measure, its universal optimality under model
(\ref{eqn:0413}) implies its universal optimality under model (\ref
{eqn:728}). \textup{(iii)} Given any criterion function satisfying conditions
\textup{(C.1)--(C.3)}, the efficiency of any measure under model (\ref
{eqn:0413}) is at least its efficiency under model (\ref{eqn:728}).
\end{theorem}

\begin{pf} (i) is readily proved by the direct comparison between
equations (\ref{eqn:03252})--(\ref{eqn:03254}) and equations (\ref
{eqn:04132})--(\ref{eqn:04133}).

For a pseudo symmetric measure, say $\xi$, it is universally optimal
for the two models as long as it maximizes the traces of the
information matrices, that is, $\operatorname{tr}(C_{\xi})=\min_{x\in\R^2}q_{\xi}(x)$
and $\operatorname{tr}(\tilde{C}_{\xi})=\min_{z\in\R}\tilde{q}_{\xi}(z)$,
respectively. If $\xi$ is also dual, $q_{\xi}(x)$ is a function
symmetric about the line of $x_1=x_2$ in view of (\ref{eqn:04243}) and
(\ref{eqn:04244}). This indicates that $\min_{x\in\R^2}q_{\xi
}(x)=\min_{z\in\R}\tilde{q}_{\xi}(z)$. Hence, the universal optimality under
the two models will be equivalent for such measure, and thus (ii) follows.

Since the information matrices of universally optimal designs are the
same for the two models in view of Proposition~\ref{prop:323} and
Theorem~\ref{thm:0413}, hence (iii) is verified as long as we can show
%
\begin{equation}
\label{eqn:04134} C_d \leq\tilde{C}_d,
\end{equation}
for any design $d$. To see (\ref{eqn:04134}), note that the column
space of $L_d+R_d$ is a subset of the column space of $[L_d|R_d]$,
hence we have $\proj(VU|VL_d|VR_d)\leq\proj(VU|V(L_d+R_d))$. Now
(\ref{eqn:04134}) follows in view of (\ref{eqn:912}) and $\tilde
{C}_d=\break T_d'V'\proj(VU| V(L_d+R_d))VT_d$.
\end{pf}

\begin{corol}\label{cor:0418}
\textup{(i)} A measure with $C_{d\xi00}$, $\tilde{C}_{\xi01}$ and $\tilde
{C}_{\xi11}$ being completely symmetric is universally optimal under
model (\ref{eqn:0413}) if and only if
%
\begin{eqnarray}
\sum_{s\in{\cal T}}p_s\frac{\partial\tilde{q}_s(z)}{\partial z}
\bigg|_{z=z^*}&=&0,\label{eqn:0424}
\\
\sum_{s\in{\cal T}}p_s&=&1.\label{eqn:04242}
\end{eqnarray}

\textup{(ii)} When $\Sigma$ is persymmetric, a pseudo symmetric dual measure is
universally optimal under model (\ref{eqn:728}) if and only if (\ref
{eqn:0424}) and (\ref{eqn:04242}) holds.
\end{corol}

\begin{remark}
Since $\tilde{q}_s(z)$ is a univariate function, one can use the
Kushner's (\citeyear{Kus97}) method to find $z^*$ and ${\cal T}$ with the
computational complexity of $O(m^2)$, where $m$ is the total number of
symmetric blocks. If we have to deal with multivariate functions such
as $q_s(x)$ (e.g., when $\Sigma$ is not persymmetric and the side
effects are directional), the computation of $x^*$ and ${\cal T}$ is
more involved but manageable. See \citet{BaiDru14} for an
example where $x$ is $5$-dimensional. Alternatively, one can build an
efficient algorithm (see the \hyperref[app]{Appendix}) based on~(\ref{eqn:3232}) to
derive the optimal measure, which further induces $x^*$ and ${\cal T}$.
\end{remark}

\section{The set ${\cal T}$ for type-$H$ covariance matrix}\label{sec5}
By restricting to the type-$H$ covariance matrix $\Sigma$, we derive
theoretical results regarding ${\cal T}$ for $2\leq t<k$. Note that the
cases of $3\leq k\leq4$ and $5\leq k\leq t$ have been studied by
\citet{KunMar00} and \citet{KunMer11}. Two special
cases of type-$H$ covariance matrix are the identity matrix and a
completely symmetric matrix.

\begin{theorem}\label{thm:0426}
Assume $\Sigma$ to be of type-$H$. \textup{(i)} If $2\leq t\leq k-2$, we have
\begin{eqnarray*}
z^*&=&0,
\\
y^*&=&k(t-1)/t-v(t-v)/k t,
\\
{\cal T}&=&\{s\dvtx f_{s,m}=u \mbox{ or }u+1,1\leq m\leq t\},
\end{eqnarray*}
where $u$ and $v$ are the integers satisfying $k=ut+v$ and $0\leq v<t$.

\textup{(ii)} If $2\leq t=k-1$, we have
%
\begin{eqnarray}
z^*&=&\frac{1}{2[k(k-3)+1/t]},\label{eqn:04245}
\\
y^*&=&k-1-\frac{2}{k}-\frac{1}{2k[k(k-3)+1/t]},
\\
{\cal T}&=&\langle s_0\rangle\cup\bigl\langle s_0'
\bigr\rangle,\label{eqn:04246}
\end{eqnarray}
where $s_0=(1\ 1\ 2\cdots t)$ and $s_0'$ is its dual sequence. Moreover, a
measure maximizes $q_{\xi}^*$ if and only if $p_{\langle s_0\rangle
}=p_{\langle s_0'\rangle}=1/2$.
\end{theorem}

\begin{pf} Due to (\ref{eqn:05233}), here we assume $\Sigma=I_k$
throughout the proof without loss of generality. For sequence
$s=(t_1\cdots t_k)$, define the quantities $\phi_s=\sum^{k-1}_{i=1}\I
_{t_i=t_{i+1}}$, $\varphi_s=\sum^{k-1}_{i=2}\I_{t_{i-1}=t_{i+1}}$,
$f_{s,m}=\sum^k_{i=1}\I_{t_i=m}$, $\chi_s=\sum^t_{m=1}f_{s,m}^2$. By
direct calculations, we have
%
\begin{eqnarray}
\tilde{q}_s(z)&=&q_{s,0}+q_{s,1}z+q_{s,2}z^2,\label{eqn:0523}
\\
q_{s,0}&=&c_{s00}=k-\chi_s/k,
\\
q_{s,1}&=&c_{s01}+c_{s02}=2(2k
\phi_s+f_{s,t_1}+f_{s,t_k}-2\chi_s)/k,
\\
q_{s,2}&=&c_{s11}+2c_{s12}+c_{s22}\label{eqn:05232}
\nonumber
\\
&=&2 \bigl[\varphi_s+k-1-(k+t-2)/k t \bigr]
\\
&&{}-2(2\chi_s-2f_{s,t_1}-2f_{s,t_k}+
\I_{t_1=t_k})/k.
\nonumber
\end{eqnarray}

(i) follows by the same approach as in Theorem~1.a of Kushner (\citeyear{Kus98})
with only more tedious arguments based on (\ref{eqn:0523})--(\ref{eqn:05232}).

Now we focus on $t=k-1$. First, we have $\phi_{s_0}=1$, $\varphi
_{s_0}=0$ and $\chi_{s_0}=k+2$, and hence $q_{s_0,0}=k-1-2/k$,
$q_{s_0,1}=-2/k$ and $q_{s_0,2}=2(k-3)+2/k t$. It can be verified that
$\tilde{q}_{s_0}(z)$ reaches its minimum at $z=z^*$. Since $\tilde
{q}_{s_0}(z)=\tilde{q}_{s_0'}(z)$, it is sufficient to show $\tilde
{q}_{s_0}(z^*)=\max_{s\in{\cal S}} \tilde{q}_s(z^*)$ for the purpose
of proving (ii).

We first restrict the consideration to the subset ${\cal S}_1=\{
s\dvtx t_1\neq t_k, s\in{\cal S}\}$. If we only exchange the treatments in
locations $\{2,\ldots,k-1\}$, the values of $\chi_s$, $f_{s,t_1}$ and
$f_{s,t_k}$ remain invariant. Note that $\tilde{q}_s(z^*)$ is
increasing in the quantity $\phi_s+2^{-1}z^*\varphi_s$. If for a
certain location, say $i$, we have $t_{i-1}=t_{i+1}\neq t_i$. At least
one of $i-1$ and $i+1$ would be in the set $\{2,\ldots,k-1\}$. After
switching this location with location $i$, $\phi_s$ will be increased
by $1$, and at the same time the amount of decrease for $\varphi_s$
will be at most $2$. Note that $z^*/2\leq1/2$ for all $p\geq3$ and
$t\geq2$, and hence a sequence, say~$s$, which maximizes $\tilde
{q}_s(z^*)$ should be of the format $s=(1_{f_{s,1}}'1|\cdots
|1_{f_{s,h}}'h)$, without loss of generality. Here, $h:=h(s)$ is the
number of distinct treatments in sequence $s$ and $\sum^h_{i=1}f_{s,i}=p$. Among sequences of this particular format, the
sequence which maximizes $\tilde{q}_s(z^*)$ should satisfy $\min
(f_{s,1},f_{s,h})\geq\max_{2\leq i\leq h-1}f_{s,i}$, where we take the
maximization over the empty set to be $0$. Without loss of generality,
we assume $t_1=\max_{1\leq i\leq t}f_{s,i}$. Now we shall show
$f_{s,1}\leq2$ for maximizing sequences as follows. Suppose
$f_{s,1}\geq3$, this indicates $h<t$. By decreasing $f_{s,1}$ by one
and changing $f_{s,h+1}$ from $0$ to $1$, the quantity $\tilde
{q}_{s}(z^*)$ is increased by the amount of
\[
\Delta_s=\frac{2}{k} \bigl[f_{s,1}-1+(4f_{s,1}-5-2k)z^*+(4f_{s,1}-8-k)
\bigl(z^*\bigr)^2 \bigr].
\]
If $k=3$, we have $\Delta_s>0$ in view of $z^*>0$ and $f_{s,1}\geq3$.
Suppose $k\geq4$, we have $0<z^*\leq(2k)^{-1}$, hence we have
\begin{eqnarray*}
k\Delta _s/2&=&f_{s,1}-1-2kz^*-p\bigl(z^*
\bigr)^2+(4f_{s,1}-5)z^*+(4f_{s,1}-8) \bigl(z^*
\bigr)^2
\\
&>&f_{s,1}-2-(4k)^{-1}>0.
\end{eqnarray*}
At this point, we have shown $\tilde{q}_{s_0}(z^*)=\max_{s\in{\cal
S}_1} \tilde{q}_s(z^*)$. By similar arguments, one can show that the
sequence $s_1=(1\ 2\cdots t1)$ maximizes $\tilde{q}_s(z^*)$ among
$s\notin{\cal S}_1$. By direct calculations, we have
\begin{eqnarray*}
\tilde{q}_{s_0}\bigl(z^*\bigr)-\tilde{q}_{s_1}\bigl(z^*
\bigr)&=&(4-2/k)z^*-4\bigl(z^*\bigr)^2/k
\\
&\geq&z^*\bigl(10/3-2/k^2\bigr)>0.
\end{eqnarray*}
Hence, (\ref{eqn:04245})--(\ref{eqn:04246}) are proved. For the rest of
(ii), the sufficiency of $p_{\langle s_0\rangle}=p_{\langle
s_0'\rangle}=1/2$ is indicated by the proof of Theorem~\ref{thm:0424}.
For the necessity, it is enough to note that the two components of
$\nabla q_{\xi}(x^*)=2(\ell_{\xi}+Q_{\xi}x^*)=2\sum_{s\in{\cal
S}}p_s(\ell_{s}+Q_{s}x^*)$ will not be identical if $p_{\langle
s_0\rangle}\neq p_{\langle s_0'\rangle}$. Hence, the lemma is concluded.
\end{pf}

\section{Examples}\label{sec6}
\begin{figure}

\includegraphics{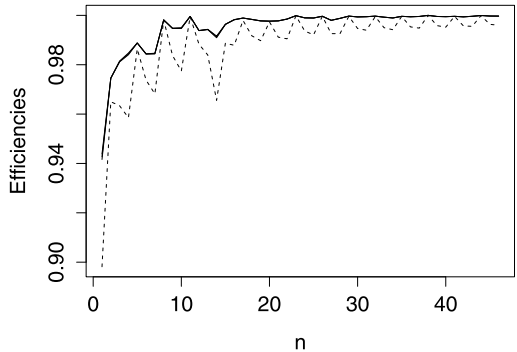}

\caption{The efficiencies of exact designs for $5\leq n\leq50$ when
$k=4$, $t=3$ and $\eta=0.5$. The $E$-efficiency is plotted by the dashed
line, while $A$-, $D$- and $T$-efficiencies are all plotted by the same solid
line.}\label{fig1}\vspace*{-7pt}
\end{figure}

This section tries to illustrate the theorems of this paper through
several examples for various combinations of $k,t,n$ and $\Sigma$. By
Theorem~\ref{thm:0424}(iii), the efficiency of a design is higher
under model (\ref{eqn:0413}) than under model (\ref{eqn:728}) for any
criterion function $\Phi$ satisfying (C.1)--(C.3) under a mild
condition, that is, $\Sigma$ is persymmetric. Hence, it is sufficient
to propose optimal or efficient designs under model (\ref{eqn:728}).
The existence of the universally optimal measure in ${\cal P}$ is
obvious in view of Lemmas \ref{lemma:0223} and \ref{lemma:0205}.
However, to derive an exact design, one has to restrict the
consideration to the subset ${\cal P}_n=\{\xi\in{\cal P}\dvtx n\xi\mbox{ is
a vector of integers}\}$. Universally, optimal measure does not
necessarily exist in ${\cal P}_n$ except for certain combinations of
$k,t,n$. In this case, one can convert $p_s$ in the equations of
Theorem~\ref{thm:0325} into $n_s$ by multiplying both sides of the
equations by $n$. Then one can define a distance between two sides of
the equations and find the solution, say $\{n_s,s\in{\cal T}\}$, to
minimize this distance. If there is universally optimal measure in
${\cal P}_n$, such approach automatically locates the universally
optimal exact design; otherwise, the exact designs thus found are
typically highly efficient under the different criteria.
See \citet{Zhe13N1} and Figure~\ref{fig1} for evidence.

Let $0\leq a_1\leq a_2\leq a_{t-1}$ be the $t$ eigenvalues of $C_d$ for
an exact design $d$. If $d$ is universally optimal, we have
$a_i=ny^*/(t-1),1\leq i\leq t-1$. Here, we define $A$-, $D$-, $E$- and
$T$-efficiencies of design $d$ as follows:
\begin{eqnarray*}
{\cal E}_A(d)&=&\frac{t-1}{ny^*}\frac{t-1}{ (\sum^{t-1}_{i=1}a_i^{-1} )}=
\frac{(t-1)^2}{n y^* (\sum^{t-1}_{i=1}a_i^{-1} )},
\\
{\cal E}_D(d)&=&\frac{t-1}{ny^*} \Biggl(\prod_{i=1}^{t-1}a_i
\Biggr)^{1/(t-1)},
\\
{\cal E}_E(d)&=&\frac{(t-1)a_1}{ny^*},
\\
{\cal E}_T(d)&=&\frac{t-1}{ny^*} \Biggl(\frac{1}{t-1}\sum
^{t-1}_{i=1}a_i \Biggr)=
\frac{\sum^{t-1}_{i=1}a_i}{ny^*}.
\end{eqnarray*}
It is well known that a universally optimal measure has unity
efficiency under these four criteria.

We begin with the discussion on the case when $\Sigma$ is of type-$H$.
For the latter, \citet{KunMar00} studied the conditions on
$p_{\langle s\rangle}$ for a pseudo symmetric design to be universally
optimal for $k=3$ and $4$, which was further extended by \citet{KunMer11} to $t\geq k\geq5$. We would comment on these cases and
then explore the case of $k\geq5$ and $t<k$. Finally, irregular form
of $\Sigma$ will be briefly discussed.

For $(k,t)=(4,2)$, Corollary~\ref{cor:0418} indicates that the
necessary and sufficient condition for a pseudo symmetric design to be
universally optimal is $p_{\langle(1\ 1\ 2\ 2)\rangle}=3p_{\langle
(1\ 2\ 1\ 2)\rangle}+p_{\langle(1\ 2\ 2\ 1)\rangle}$. Theorem~2 of \citet{KunMar00} proposed $p_{\langle(1\ 1\ 2\ 2)\rangle}=p_{\langle
(1\ 2\ 2\ 1)\rangle}=1/2$, which is sufficient but not necessary for
universal optimality. For $k=3$ and $(k,t)=(4,3)$, Corollary~\ref
{cor:0418} indicates that sufficient conditions regarding $p_{\langle
s\rangle}$ given by Theorems 1 and 3 of \citet{KunMar00} are
also necessary.

For $t\geq k=4$, \citet{KunMar00} showed that the optimal
values of $p_{\langle s \rangle}$ are given by irrational numbers, and
hence an exact universally optimal design does not exist. In fact,
based on Theorem~\ref{thm:0325} here, one can derive efficient exact
designs for the majority values of $t$ and $n$. For example, $d_1$
below with $t=4$ and $n=10$ yields the efficiencies of ${\cal
E}_A(d_1)=0.9943$, ${\cal E}_D(d_1)=0.9946$, ${\cal E}_E(d_1)=0.9682$
and ${\cal E}_T(d_1)=0.9949$. Note that the $E$-efficiency is relatively
lower than other efficiencies due to the asymmetry of the design.
\[
d_1=\left[ %
\matrix{ 2& 1 &
4& 3 & 1& 1 & 3& 2 & 4& 3
\cr
2 & 1 & 4 & 3 & 2 & 4 & 4 & 3 & 2 & 2
\cr
1 & 3 & 3 & 1 & 4 & 3 & 2 & 1 & 1 & 4
\cr
1 & 4& 2 & 2& 4 & 3& 2 & 4& 3 & 1 }
 \right].
\]

For $t\geq k\geq5$, \citet{KunMer11} showed that the set
${\cal T}$ should include sequences $(1\ 2\cdots k)$, $(1\ 1\ 2\cdots
k-3\ k-2\ k-2)$, $s_0$ and its dual sequence $s_0'$ as defined in
Theorem~\ref{thm:0426}. The optimal proportion for them are again
irrational numbers. Further, they proposed the use of type I orthogonal
array ($\mathit{OA}_I$), that is, $p_{\langle(1\ 2\cdots k) \rangle}=1$, and proved
that the $T$-efficiencies of such designs are at least $0.94$. Note that
$\mathit{OA}_I$ is pseudo symmetric, hence its efficiencies are identical under
criteria $A$, $D$, $E$ and $T$.

When $t=k-1$, Theorem~\ref{thm:0426}(ii) indicates that a pseudo
symmetric design with $p_{\langle s_0\rangle}=p_{\langle s_0'\rangle
}=1/2$ will be universally optimal. For example, when $t=4$ and $k=5$,
$d_2$ below with $n=24$ is universally optimal. Here, the first $12$
sequences are equivalent to $(1\ 1\ 2\ 3\ 4)$ while the rest are equivalent to
$(1\ 2\ 3\ 4\ 4)$.
\begin{eqnarray*}
d_2&=&\left[ \matrix{
1&1&1&2&2&2&3&3&3&4&4&4
\cr
1&1&1&2&2&2&3&3&3&4&4&4
\cr
4&2&3&1&4&3&1&4&2&1&2&3
\cr
2&3&4&4&3&1&2&1&4&3&1&2
\cr
3&4&2&3&1&4&4&2&1&2&3&1 }
 \right.\\
 &&\hspace*{5pt}\left.\matrix{
3&4&2&3&1&4&4&2&1&2&3&1\cr
2&3&4&4&3&1&2&1&4&3&1&2\cr
4&2&3&1&4&3&1&4&2&1&2&3\cr
1&1&1&2&2&2&3&3&3&4&4&4\cr
1&1&1&2&2&2&3&3&3&4&4&4 }\right].
\end{eqnarray*}

When $2\leq t<k-1$, there is a large variety of symmetric blocks in
${\cal T}$ and there will be infinity many solutions for optimal
sequence proportions. Even for $t=2$ and $k=5$, we shall have ${\cal
T}=\langle(1\ 1\ 1\ 2\ 2) \rangle\cup\langle(1\ 1\ 2\ 2\ 2)\rangle\cup\langle
(1\ 1\ 2\ 1\ 2)\rangle\cup\langle(1\ 2\ 1\ 2\ 2)\rangle\cup\langle(1\ 1\ 2\ 2\ 1)\rangle
\cup
\langle(1\ 2\ 2\ 1\ 1)\rangle\cup\langle(1\ 2\ 1\ 1\ 2)\rangle\cup\langle
(1\ 2\ 2\ 1\ 2)\rangle\cup\langle(1\ 2\ 1\ 2\ 1)\rangle\cup\langle(1\ 2\ 2\ 2\ 1)\rangle$.
Let $p_1,\ldots,p_{10}$ be the proportions of these symmetric blocks. A
pseudo symmetric design with $p_1=p_2$, $p_3=p_4$, $p_5=p_6$,
$p_7=p_8$, $1.8(p_1+p_2)=2.2(p_3+p_4+p_7+p_8)+4p_9+0.4p_{10}$, $\sum^{10}_{i=1}p_i=1$ and $p_i\geq0$ will be universally optimal. One
simple solution is $p_5=p_6=1/2$. Hence a design which assigns $1/4$ of
its blocks to sequences $(1\ 1\ 2\ 2\ 1)$, $(2\ 2\ 1\ 1\ 2)$, $(1\ 2\ 2\ 1\ 1)$ and $(2\ 1\ 1\ 2\ 2)$
is universally optimal.

At last, we would like to convey the message that the deviation of
$\Sigma$ from type-$H$ has large impact on the choice of designs. For
simplicity of illustration, we consider the form $\Sigma=(\I
_{i=j}+\eta
\I_{|i-j|=1})_{1\leq i,j\leq k}$. When $k=t=5$ and $\eta=0.5$, the
efficiency of $\mathit{OA}_I$ reduces to $0.8232$. In fact, Corollary~\ref
{cor:0418} indicates that $\langle(1\ 1\ 2\ 3\ 3)\rangle$, instead of
$\langle
(1\ 2\ 3\ 4\ 5)\rangle$ for $\eta=0$, becomes the dominating symmetric block
among the four. To be more specific, a pseudo symmetric design with
sequences solely from $\langle(1\ 1\ 2\ 3\ 3)\rangle$ yields the efficiency of
$0.9999$ for all four criteria. When we tune $\eta$ to $0.9$, the
efficiency of $\mathit{OA}_I$ further reduces to $0.3395$, while the symmetric
design based on $\langle(1\ 1\ 2\ 3\ 3)\rangle$ becomes even more efficient.
One the other hand, when $\eta$ takes negative values, the efficiency
of $\mathit{OA}_I$ becomes even higher than $0.94$. Similar phenomena are
observed for other cases of $t\geq k$.

For $t<k$, we also observe that the value of $\eta$ influences the
choice of design substantially. The details are omitted due to the
limit of space. We end this section by Figure~\ref{fig1}. It shows that
the linear equations system in Theorem~\ref{thm:0325} is powerful in
deriving efficient exact designs for arbitrary values of $n$.

\begin{appendix}\label{app}
\section*{Appendix: The algorithm based on Theorem~\lowercase{\protect\ref{thm:323}}}
Recall that $m$ is be total number of distinct symmetric blocks and
$s_1,s_2,\break \ldots,s_m$ are the $m$ representatives\vadjust{\goodbreak} for each of the symmetric
blocks. Note that two pseudo symmetric measures with the same vector of
$P_{\xi}=(p_{\langle s_1\rangle},p_{\langle s_2\rangle
},\ldots,\break p_{\langle
s_m\rangle})$ have the same information matrix and hence the same
performance under all optimality criteria. For a measure $\xi$ and a
sequence $s$, we define
%
\begin{equation}
\label{eqn:409} \theta(P_{\xi},s)=\operatorname{tr}\bigl(R_sR_{\xi}^{-1}
\bigr)-\operatorname{tr}\bigl(Q_sQ_{\xi}^{-1}\bigr).
\end{equation}
We also define $\theta^*(P_{\xi})=\max_{1\leq i\leq m}\theta
(P_{\langle
d\rangle},s_i)$ and $e_i$ to be vector of length $m$ with the $i$th
entry as 1 and other entries as 0.

\textit{Step} 0: Choose tuning parameters $\epsilon>0$ and $\omega$ such that
$\epsilon$ is in a small neighborhood of zero and $\omega$ is in a
neighborhood of one.

\textit{Step} 1: Choose initial measure $P^{(0)}=P_{\xi_0}$. Put $i_0=\argmin
_{1\leq j\leq m}\theta^*(e_i)$ and $n=0$, then let $P^{(0)}=e_{i_0}$.

\textit{Step} 2: Check optimality. If $\theta_n:=\theta^*(P^{(n)})>1+\epsilon$,
go to step 3. Otherwise, output the optimal measure as $P^{(n)}$.

\textit{Step} 3: Update the measure. Let $i_{n+1}=\argmax_{1\leq i\leq m}\theta
(P^{(n)},s_i)$ and the updated measure is $P^{(n+1)}=(\theta
_n-1)^{\omega}e_{i_{n+1}}+(1-(\theta_n-1)^{\omega})P^{(n)}$. Increase
$n$ by 1 and go back to step 2.

\begin{remark}
There is a possibility of tie in choosing $i_0$ in step 1 and $i_{n+1}$
in step 3. The strategy in such case is quite arbitrary. Let $\Xi_n=\{
i\dvtx \theta(P^{(n)},s_i)=\theta^*(P^{(n)})\}$. If $|\Xi_n|>1$, one can
either choose an arbitrary $j_n\in\Xi_n$ and let $i_{n+1}=j_n$ or
replace $e_{i_{n+1}}$ in step 3 by $|\Xi_n|^{-1}\sum_{i\in\Xi_n}e_i$.
The same strategy applies to the choice of $i_0$.
\end{remark}

\begin{remark}
Note that the update algorithm in step 3 is essentially a steepest
descent algorithm. The parameter $\omega$ is to adjust for the length
of step for the \textit{best} direction. By the concavity of the
optimality criteria, the global optimum is guaranteed to be found. In
the examples of this paper, $\omega=1$ works well enough. The parameter
$\epsilon$ is used to adjust for time of convergence. When the
sequential algorithm converges very slow, one can increase $\epsilon$
to save time. In most examples of this paper, setting $\epsilon
=10^{-7}$ enable us to obtain the optimal design within $10$ seconds.
\end{remark}
\end{appendix}

\section*{Acknowledgments}
We are grateful to the Associate Editor of this paper and two referees
for their constructive comments on earlier versions of this manuscript.




\printaddresses
\end{document}